\documentclass[11pt]{article}
\usepackage[]{amsmath,amssymb}

\newtheorem{theorem}{Theorem}[section]
\newtheorem{lemma}[theorem]{Lemma}
\newtheorem{proposition}[theorem]{Proposition}

\newtheorem{exAux}[theorem]{Example}

\newtheorem{Def}[theorem]{Definition}
\newenvironment{definition}{\begin{Def} \rm}{\end{Def}}
\newtheorem{Note}[theorem]{Note}
\newenvironment{note}{\begin{Note} \rm}{\end{Note}}
\newtheorem{Problem}[theorem]{Problem}

\newtheorem{Rem}[theorem]{Remark}

\newtheorem{Not}[theorem]{Notation}

\newtheorem{Ass}[theorem]{Assumption}

\newenvironment{proof}{\medskip\noindent{\bf Proof.\ }}{\qed\medskip}

\newcommand{\qed}{\hfill\mbox{$\Box$\qquad\qquad}}

\renewcommand{\indent}{\hspace{6mm}}

\begin{document}
\thispagestyle{empty}

\begin{center}
\LARGE \bf
\noindent
The split decomposition of a tridiagonal pair
\end{center}

\smallskip

\begin{center}
\Large
Kazumasa Nomura and Paul Terwilliger
\end{center}

\smallskip

\begin{quote}
\small 
\begin{center}
\bf Abstract
\end{center}
\indent
Let $\mathbb{K}$ denote a field and let $V$ denote a vector space
over $\mathbb{K}$ with finite positive dimension.
We consider a pair of linear transformations 
$A:V \to V$ and $A^*:V \to V$ that satisfy (i)--(iv) below:
\begin{itemize}
\item[(i)] 
Each of $A$, $A^*$ is diagonalizable.
\item[(ii)] 
There exists an ordering $V_{0},V_{1},\ldots,V_{d}$ 
of the eigenspaces of $A$ such that
$A^* V_i \subseteq V_{i-1} + V_{i} + V_{i+1}$ for $0 \leq i \leq d$,
where $V_{-1}=0$, $V_{d+1}=0$.
\item[(iii)]
There exists an ordering $V^*_{0},V^*_{1},\ldots,V^*_{\delta}$
of the eigenspaces of $A^*$ such that
$A V^*_i \subseteq V^*_{i-1} + V^*_{i} + V^*_{i+1}$ for $0 \leq i \leq \delta$,
where $V^*_{-1}=0$, $V^*_{\delta+1}=0$.
\item[(iv)]
There is no subspace $W$ of $V$ such that both
$AW \subseteq W$, $A^* W \subseteq W$, other than $W=0$ and $W=V$.
\end{itemize}
We call such a pair a {\em tridiagonal pair} on $V$.
In this note we obtain two results. 
First, we show that each of $A,A^*$ is determined
up to affine transformation by the $V_i$ and $V^*_i$.
Secondly,
we characterize the case in which the $V_i$ and $V^*_i$
all have dimension one.  
We prove both results using a certain decomposition of
$V$ called the split decomposition.
\end{quote}

\section{Introduction}

\indent
Throughout this note $\mathbb{K}$ will denote a field
and $V$ will denote a vector space over $\mathbb{K}$ with finite
positive dimension.
Let $\text{End}(V)$ denote the $\mathbb{K}$-algebra of all 
$\mathbb{K}$-linear transformations from $V$ to $V$.

\medskip

For $A \in \text{End}(V)$ and for a subspace $W \subseteq V$,
we call $W$ an {\em eigenspace} of $A$ 
whenever $W \neq 0$ and there exists $\theta \in \mathbb{K}$
such that $W=\{v \in V \,|\, Av=\theta v\}$.
We say $A$ is {\em diagonalizable} whenever $V$ is spanned by
the eigenspaces of $A$.
We now recall the notion of a tridiagonal pair.

\begin{definition} \cite{ITT}     \label{def:TDpair}   \samepage
By a {\em tridiagonal pair} on $V$ we mean an ordered pair
of elements $A,A^*$ taken from $\text{End}(V)$ that satisfy (i)--(iv) below:
\begin{itemize}
\item[(i)] 
Each of $A$, $A^*$ is diagonalizable.
\item[(ii)] 
There exists an ordering $V_{0},V_{1},\ldots,V_{d}$ 
of the eigenspaces of $A$ such that
\[
  A^* V_i \subseteq V_{i-1} + V_{i} + V_{i+1}  \qquad\qquad (0 \leq i \leq d),
\]
 where $V_{-1}=0$, $V_{d+1}=0$.
\item[(iii)] 
There exists an ordering $V^*_{0},V^*_{1},\ldots,V^*_{\delta}$
of the eigenspaces of $A^*$ such that
\[
A V^*_i \subseteq V^*_{i-1} + V^*_{i} + V^*_{i+1} \qquad\qquad
            (0 \leq i \leq \delta),
\]
where $V^*_{-1}=0$, $V^*_{\delta+1}=0$.
\item[(iv)] 
There is no subspace $W$ of $V$ such that both
$AW \subseteq W$, $A^* W \subseteq W$, other than $W=0$ and $W=V$.
\end{itemize}
\end{definition}

\begin{note}              \samepage
It is a common notational convention to use $A^*$ to represent the
conjugate-transpose of $A$. We are not using this convention.
In a tridiagonal pair $A,A^*$ the linear transformations $A$ and
$A^*$ are arbitrary subject to (i)--(iv) above.
\end{note}

\medskip

We refer the reader to
\cite{ITT}, \cite{IT:shape},
\cite{N:refine}
for background on tridiagonal pairs.

\medskip

Referring to Definition \ref{def:TDpair} we have 
$d=\delta$ \cite[Lemma 4.5]{ITT}; 
we call this common value  the {\em diameter} of $A,A^*$.
For $0 \leq i \leq d$ the dimensions of $V_i$ and $V^*_i$ coincide;
we denote this common value by $\rho_i$, and observe that $\rho_i \neq 0$.
The sequence $\rho_0,\rho_1,\ldots,\rho_d$ is symmetric and
unimodal; i.e. $\rho_i=\rho_{d-i}$ for $0 \leq i \leq d$ 
and $\rho_{i-1} \leq \rho_i$ for $1 \leq i \leq d/2$
\cite[Corollaries 5.7, 6.6]{ITT}.
We call the vector $(\rho_0,\rho_1, \ldots,\rho_d)$ 
the {\em shape} of $A,A^*$.
By a {\em Leonard pair} we mean a tridiagonal pair with shape 
$(1,1,\ldots,1)$ \cite[Definition 1.1]{T:Leonard}.
See
\cite{T:Leonard},
\cite{T:qRacah},
\cite{T:survey}
for background information on Leonard pairs.

\medskip

In this note we obtain the following two results. 
Let $A,A^*$ denote a tridiagonal pair from Definition \ref{def:TDpair}.
First, we show that each of $A,A^*$ is determined
up to affine transformation by the $V_i$ and $V^*_i$.
Secondly, we characterize the Leonard pairs among the
tridiagonal pairs. We prove both results using a certain
decomposition of the underlying vector space called 
the split decomposition \cite[Section 4]{ITT}.

\section{The split decomposition}

\indent
In this section we recall the split decomposition \cite[Section 4]{ITT}.
We start with a comment. 
Referring to Definition \ref{def:TDpair}, 
since $V_0,V_1, \ldots, V_d$ are the eigenspaces of $A$ and since
$A$ is diagonalizable we have
\begin{equation}        \label{eq:sumVi}
  V = V_0+V_1+\cdots + V_d   \qquad\qquad   \text{(direct sum)}.
\end{equation}
Similarly
\begin{equation}        \label{eq:sumVsi}
  V = V^*_0+V^*_1+\cdots + V^*_d \qquad\qquad  \text{(direct sum)}. 
\end{equation}
For $0 \leq i \leq d$ define
\begin{equation}                   \label{eq:defUi}
 U_i=(V^*_0+V^*_1+\cdots+V^*_i) \cap (V_i+V_{i+1}+\cdots+V_d).
\end{equation}
By \cite[Theorem 4.6]{ITT},
\begin{equation}               \label{eq:decomp}
   V=U_0+U_1+\cdots+U_d   \qquad\qquad \text{(direct sum)},
\end{equation}
and for $0 \leq i \leq d$ both
\begin{eqnarray}
 U_0+U_1+\cdots+U_i &=& V^*_0+ V^*_1+\cdots+V^*_i,    \label{eq:sumUiVsi} \\
 U_i+U_{i+1}+\cdots+U_d &=& V_i+V_{i+1}+\cdots+V_d.   \label{eq:sumUiVi}
\end{eqnarray}
For $0 \leq i \leq d$ let $\theta_i$ (resp. $\theta^*_i$)
denote the eigenvalue of $A$ (resp. $A^*$) associated with
the eigenspace $V_i$ (resp. $V^*_i$).
Then by \cite[Theorem 4.6]{ITT} both
\begin{eqnarray}
  (A-\theta_i I)U_i &\subseteq& U_{i+1},          \label{eq:AUi}  \\
  (A^*-\theta^*_i I)U_i &\subseteq& U_{i-1},      \label{eq:AsUi}
\end{eqnarray}
where $U_{-1}=0$ and $U_{d+1}=0$.
The sequence 
$U_0,U_1,\ldots, U_d$ is called the {\it split decomposition} of $V$
\cite[Section 4]{ITT}.

\section{A subalgebra of $\text{End}(V)$}

\indent
The following subalgebra of $\text{End}(V)$ will be useful to us.
Referring to Definition \ref{def:TDpair},
let $\cal D$ denote the subalgebra of $\text{End}(V)$ generated by $A$. 
In what follows we often view $\cal D$ as a vector space over $\mathbb{K}$.
The dimension of this vector space is $d+1$ since $A$ is diagonalizable 
with $d+1$ eigenspaces.
Therefore $\{A^i \,|\, 0 \leq i \leq d\}$ is a basis for $\cal D$. 
There is another basis for $\cal D$ that is better suited to our purpose. 
To define it we use the following notation.
Let $\mathbb{K}[\lambda]$ denote the $\mathbb{K}$-algebra of all polynomials 
in an indeterminate $\lambda$ that have coefficients in $\mathbb{K}$. 
For $0 \leq i \leq d$ we define $\tau_i \in \mathbb{K}[\lambda]$ by
\begin{equation}              \label{eq:deftau}
  \tau_i = (\lambda - \theta_0)(\lambda-\theta_1)\cdots(\lambda-\theta_{i-1}).
\end{equation}
We note that $\tau_i$ is monic with degree $i$.
Therefore $\{\tau_i(A) \,|\, 0 \leq i \leq d \}$ 
is a basis for $\cal D$.
Combining (\ref{eq:AUi}) and (\ref{eq:deftau}) we find
\begin{equation}                      \label{eq:tauiAU0}
   \tau_i(A)U_0 \subseteq U_i \qquad\qquad (0 \leq i \leq d).
\end{equation}

\medskip

The following lemma is a variation on  \cite[Lemma 6.5]{ITT}; 
we give a short proof for the convenience of the reader.

\medskip

\begin{lemma}        \label{lem:calDu}    \samepage 
Referring to Definition \ref{def:TDpair}, 
for all nonzero $u \in V^*_0$ and nonzero $X \in {\cal D}$,
we have $Xu \neq 0$.
\end{lemma}

\begin{proof}
It suffices to show that the vector spaces
$\cal D$ and ${\cal D}u$ have the same dimension.
We saw earlier that 
$\{\tau_i(A) \,|\, 0 \leq i \leq d\}$ is a basis for $\cal D$.
We show that 
$\{\tau_i(A)u \,|\, 0 \leq i \leq d\}$ is a basis for ${\cal D}u$.
By (\ref{eq:decomp}), (\ref{eq:tauiAU0}),  and since $U_0=V^*_0$, 
this will hold if we can show
$\tau_i(A)u \neq 0$ for $0 \leq i \leq d$.
Let $i$ be given and suppose $\tau_i(A)u=0$. We will obtain
a contradiction by displaying a subspace $W$ of $V$ that violates
Definition \ref{def:TDpair}(iv). 
Observe that $i \neq 0$ since $\tau_0=1$ and $u \neq 0$; therefore $i\geq 1$.
By (\ref{eq:deftau}) and since $\tau_i(A)u=0$ we find
$u \in V_0 + V_1+\cdots + V_{i-1}$, so
\begin{equation}           \label{eq:lemaux1}
   u \in V^*_0 \cap (V_0+V_1+\cdots+V_{i-1}).
\end{equation}
Define 
\begin{equation}           \label{eq:lemaux2}
  W_r = (V^*_0+V^*_1+\cdots+V^*_r) \cap (V_0+V_1+\cdots+V_{i-r-1})
\end{equation}
for $0 \leq r \leq i-1$ and put
\begin{equation}              \label{eq:sumW}
  W = W_0+W_1+\cdots + W_{i-1}.
\end{equation}
We show $W$ violates Definition \ref{def:TDpair}(iv).
Observe that $W \neq 0$ since the nonzero vector $u \in W_0$ by 
(\ref{eq:lemaux1}) and since $W_0 \subseteq W$.
Next we show $W \neq V$.
By (\ref{eq:lemaux2}), for $0 \leq r \leq i-1$  we have
\begin{eqnarray*}
 W_r &\subseteq& V^*_0 + V^*_1 + \cdots +V^*_r  \\
     &\subseteq& V^*_0 + V^*_1 + \cdots +V^*_{i-1}.
\end{eqnarray*}
By this and (\ref{eq:sumW}) we find
\begin{eqnarray*}
 W &\subseteq&  V^*_0 + V^*_1+ \cdots +V^*_{i-1}  \\
   &\subseteq&  V^*_0 + V^*_1+ \cdots +V^*_{d-1}.
\end{eqnarray*}
Combining this with (\ref{eq:sumVsi}) and using $V^*_d \neq 0$ we find
$W \neq V$. 
We now show $AW \subseteq W$.
To this end, we show that 
$(A-\theta_{i-r-1}I)W_r \subseteq W_{r+1}$ for $0 \leq r \leq i-1$, 
where $W_i=0$.
Let $r$ be given.
From the construction we have
\begin{equation}                 \label{eq:L1}
   (A - \theta_{i-r-1} I) \sum_{h=0}^{i-r-1}V_h 
    = \sum_{h=0}^{i-r-2} V_h.
\end{equation}
By Definition \ref{def:TDpair}(iii) we have
\begin{equation}                  \label{eq:L2}
    (A - \theta_{i-r-1} I) \sum_{h=0}^{r} V^*_h \subseteq \sum_{h=0}^{r+1}V^*_h.
\end{equation}
Combining (\ref{eq:L1}) and (\ref{eq:L2})
we find $(A-\theta_{i-r-1}I)W_{r} \subseteq W_{r+1}$ as desired.
We have shown $AW \subseteq W$.
We now show $A^*W \subseteq W$.
To this end, we show that 
$(A^*-\theta^*_r I)W_r \subseteq W_{r-1}$ for $0 \leq r \leq i-1$, 
where $W_{-1}=0$.
Let $r$ be given.
From the construction we have
\begin{equation}                   \label{eq:L3}
   (A^* - \theta^*_r I) \sum_{h=0}^{r} V^*_h 
    = \sum_{h=0}^{r-1} V^*_h.
\end{equation}
By Definition \ref{def:TDpair}(ii) we have
\begin{equation}                    \label{eq:L4}
  (A^*-\theta^*_r I) \sum_{h=0}^{i-r-1} V_h \subseteq \sum_{h=0}^{i-r}V_h.
\end{equation}
Combining (\ref{eq:L3}) and (\ref{eq:L4})
we find $(A^*-\theta^*_r I)W_{r} \subseteq W_{r-1}$ as desired.
We have shown $A^*W \subseteq W$.
We have now shown that $W \neq 0$, $W \neq V$, $AW \subseteq W$,
$A^* W \subseteq W$, contradicting Definition \ref{def:TDpair}(iv).
We conclude $\tau_i(A)u\not=0$ and the result follows.
\end{proof}

\section{Each of $A,A^*$ is determined by the eigenspaces}

\indent
Let the tridiagonal pair $A,A^*$ be as in Definition \ref{def:TDpair}.
In this section we show that each of $A,A^*$ is determined
up to affine transformation by the eigenspaces $V_i$, $V^*_i$.
Our main result is based on the following proposition.

\medskip

\begin{proposition}        \label{prop:main}   \samepage
Referring to Definition \ref{def:TDpair}, assume $d \geq 1$.
Then the following (i), (ii) are equivalent for all $X \in \text{\rm End}(V)$.
\begin{itemize}
\item[(i)] $X \in {\cal D}$ and $XV^*_0 \subseteq V^*_0+V^*_1$.
\item[(ii)] There exist scalars $r,s$ in $\mathbb{K}$ such that $X=rA+sI$.
\end{itemize}
\end{proposition}

\begin{proof}
(i)$\Rightarrow$(ii):
Assume $X \neq 0$; otherwise the result is trivial.
Pick a nonzero $u \in V^*_0$ and note that $u \in U_0$ by (\ref{eq:sumUiVsi}).
We have $Xu \in V^*_0+V^*_1$ by assumption so
\begin{equation}            \label{eq:uinU0U1}
               Xu \in U_0+U_1
\end{equation}
in view of (\ref{eq:sumUiVsi}). 
Recall 
$\{\tau_i(A) \,|\, 0 \leq i \leq d\}$
is a basis for $\cal D$. We assume $X \in {\cal D}$ so there exist
$\alpha_i \in \mathbb{K}$ $(0 \leq i \leq d)$ such that
\begin{equation}         \label{eq:B}
  X=\sum_{i=0}^d \alpha_i \tau_i(A).
\end{equation}
We show $\alpha_i=0$ for $2\leq i \leq d$.
Suppose not and define 
$\eta= \max\{i \,|\, 2 \leq i \leq d, \; \alpha_i \not=0 \}$.
We will obtain a contradition by showing
\begin{equation}         \label{eq:propaux1}
    0 \not= U_{\eta} \cap (U_0+U_1 + \cdots + U_{\eta-1}).
\end{equation}
Note that $\tau_{\eta}(A)u \not=0$ by Lemma \ref{lem:calDu} and
$\tau_{\eta}(A)u \in  U_{\eta}$ by (\ref{eq:tauiAU0}).
Also by (\ref{eq:B}) we find $\tau_{\eta}(A)u$ is in the span of
$Xu$ and $\tau_0(A)u, \tau_1(A)u,\ldots, \tau_{\eta-1}(A)u$;
combining this with (\ref{eq:tauiAU0}) and (\ref{eq:uinU0U1}) we find
$\tau_{\eta}(A)u$ is contained in 
$U_0+U_1+\cdots + U_{\eta-1}$.
By these comments $\tau_{\eta}(A)u$ is a nonzero element in 
$U_{\eta} \cap (U_0+U_1+\cdots + U_{\eta-1})$ and (\ref{eq:propaux1}) follows.
Line (\ref{eq:propaux1}) contradicts (\ref{eq:decomp}) and we conclude
$\alpha_i= 0$ for  $2 \leq i \leq d$.
Now $X = \alpha_1 \tau_1(A) + \alpha_0 I $.
Therefore $X=rA+sI$ with $r=\alpha_1$ and $s=\alpha_0-\alpha_1\theta_0$.

(ii)$\Rightarrow$(i):
Immediate from Definition \ref{def:TDpair}(iii).
\end{proof}

\medskip

The following is our first main theorem.

\medskip

\begin{theorem}              \label{thm:main}        \samepage
Let $A,A^*$ denote a tridiagonal pair on $V$, with
eigenspaces $V_i,V^*_i$ $(0 \leq i \leq d)$ as in Definition \ref{def:TDpair}.
Let $A',A^{*\prime}$ denote a second tridiagonal pair on $V$, 
with eigenspaces $V'_i, V^{*\prime}_i$  $(0 \leq i \leq d)$
as in Definition \ref{def:TDpair}.
Assume $V_i=V'_i$ and $V^*_i=V^{*\prime}_i$ for $0 \leq i \leq d$.
Then $\text{\rm Span}\{A,I\} = \text{\rm Span}\{A',I\}$ and
$\text{\rm Span}\{A^*,I\}= \text{\rm Span}\{A^{*\prime},I\}$.
\end{theorem}

\begin{proof}
Assume $d\geq 1$; otherwise the result is clear.
Let $\cal D$ (resp. ${\cal D}'$) denote the subalgebra of $\text{End}(V)$
generated by $A$ (resp. $A'$). Since $V_i=V'_i$ for $0 \leq i \leq d$
we find ${\cal D}={\cal D}'$ so $A' \in {\cal D}$.
Applying Proposition \ref{prop:main} to the tridiagonal pair $A,A^*$ 
(with $X=A'$), there exist $r, s \in \mathbb{K}$ such that $A'=rA+sI$. 
Note that $r\not=0$; otherwise $A'=sI$ has a single eigenspace
which contradicts $d\geq 1$.
It follows that $\text{Span}\{A,I\}=\text{Span}\{A',I\}$.
Similarly we find $\text{Span}\{A^*,I\}=\text{Span}\{A^{*\prime},I\}$.
\end{proof}

\section{A characterization of the Leonard pairs}

\indent
In this section we obtain a characterization of the Leonard pairs
among the tridiagonal pairs.
This characterization is based on the notion of the 
{\em switching element} of a Leonard pair \cite{NT:switch}.
We briefly recall this notion. Let the tridiagonal pair
$A,A^*$ be as in Definition \ref{def:TDpair}, and assume the corresponding
shape is $(1,1,\ldots,1)$ so that $A,A^*$ is a Leonard pair.
For $0 \leq i \leq d$ let $E_i$ denote the element of $\text{End}(V)$
such that 
$(E_i-\delta_{ij}I)V_j=0$ for $0 \leq j \leq d$.
We observe
(i) $E_iE_j = \delta_{ij}E_i$ $(0 \leq i,j\leq d)$;
(ii) $\sum_{i=0}^d E_i=I$;
(iii) $A = \sum_{i=0}^d \theta_i E_i$. 
We further observe that $E_0,E_1, \ldots, E_d$ is a basis for $\cal D$.
We define
\[
 S = \sum_{r=0}^d
      \frac{\phi_d\phi_{d-1}\cdots\phi_{d-r+1}}
           {\varphi_1\varphi_2\cdots\varphi_r} E_r,
\]
where $\varphi_1, \varphi_2, \ldots, \varphi_d$
(resp. $\phi_1, \phi_2, \ldots, \phi_d$)
is the first split sequence (resp. second split sequence)
for $A,A^*$ \cite[Definitions 3.10, 3.12] {T:Leonard}. 
The element $S$ is called the {\em switching element} for 
$A,A^*$ \cite[Definition 5.1]{NT:switch}.

\medskip

The switching element has the following property.

\medskip

\begin{theorem}    \cite[Theorem 6.7]{NT:switch}  \label{thm:S} \samepage
Let $A,A^*$ denote a Leonard pair on $V$ and let $S$ denote the
switching element for $A,A^*$.
Let $V_i,V^*_i$ be as in Definition \ref{def:TDpair}, and
let $\cal D$ denote the subalgebra of $\text{\rm End}(V)$ generated by $A$.
Then for $X \in \text{\rm End}(V)$ the following (i), (ii) are
equivalent.
\begin{itemize}
\item[(i)]  $X$ is a scalar multiple of $S$.
\item[(ii)] $X \in {\cal D}$ and $XV^*_0 \subseteq V^*_d$.
\end{itemize}
\end{theorem}

\medskip

The following is our second main result.

\medskip

\begin{theorem}         \label{thm:main2}  \samepage
Referring to Definition \ref{def:TDpair},
let $\cal D$ denote the subalgebra of $\text{\rm End}(V)$ generated by $A$.
Then the following (i), (ii) are equivalent.
\begin{itemize}
\item[(i)] There exists a nonzero $X \in {\cal D}$ 
such that $XV^*_0 \subseteq V^*_d$.
\item[(ii)] The pair $A,A^*$ is a Leonard pair.
\end{itemize}
\end{theorem}

\begin{proof}
(i)$\Rightarrow$(ii):
Recall $\{\tau_i(A) \,|\, 0 \leq i \leq d\}$ is a basis 
for $\cal D$.
Therefore there exist scalars 
$\alpha_0, \alpha_1, \ldots, \alpha_d$ in $\mathbb{K}$, not all zero,
such that
\begin{equation}      \label{eq:X}
   X= \sum_{i=0}^d \alpha_i \tau_i(A).
\end{equation}
Fix a nonzero $u \in V^*_0$.
Then $Xu \in V^*_d$ so $(A^*-\theta^*_dI)Xu=0$.
In this equation we evaluate $X$ using (\ref{eq:X}), rearrange terms,
and use $A^*u=\theta^*_0u$ to find
\begin{eqnarray*}
0 &=& (A^* - \theta^*_d I)Xu \\
 &=& \sum_{i=0}^d \alpha_i (A^* - \theta^*_d I)\tau_i(A) u    \\
 &=& \sum_{i=0}^d \alpha_i(\theta^*_i - \theta^*_d)\tau_i(A)u  
     + \sum_{i=0}^d \alpha_i(A^* - \theta^*_i I)\tau_i(A)u       \\
 &=& \sum_{i=0}^{d-1} \alpha_i(\theta^*_i - \theta^*_d)\tau_i(A)u  
     + \sum_{i=1}^d \alpha_i(A^* - \theta^*_i I)\tau_i(A)u    \\
 &=& \sum_{i=0}^{d-1} (\alpha_i(\theta^*_i-\theta^*_d)\tau_i(A)u+
                  \alpha_{i+1}(A^*-\theta^*_{i+1}I)\tau_{i+1}(A)u).
\end{eqnarray*}
In the above line, for $0 \leq i \leq d-1$ the summand at $i$ is contained in
$U_i$ in view of  (\ref{eq:AsUi}) and (\ref{eq:tauiAU0}),
so this summand is $0$ in view of (\ref{eq:decomp}).
Therefore
\begin{equation}             \label{eq:aux2}
 \alpha_i(\theta^*_i-\theta^*_d)\tau_i(A)u
           + \alpha_{i+1}(A^* - \theta^*_{i+1}I)\tau_{i+1}(A)u=0
  \qquad\qquad (0 \leq i \leq d-1).
\end{equation}
Suppose for the moment that there exists an integer $i$ $(0 \leq i \leq d-1)$ 
such that $\alpha_{i+1}=0$.
Then $\alpha_i(\theta^*_i-\theta^*_d)\tau_i(A)u=0$ by (\ref{eq:aux2}).
But $\tau_i(A)u \not=0$ by Lemma \ref{lem:calDu}
and $\theta^*_i\not= \theta^*_d$ since 
$\theta^*_0,\theta^*_1,\ldots, \theta^*_d$ are distinct, so $\alpha_i=0$.
Therefore $\alpha_{i+1}=0$ implies $\alpha_i=0$ for $0 \leq i \leq d-1$. 
By this and since 
$\alpha_0, \alpha_1, \ldots, \alpha_d$ are not all zero,
there exists an integer $j$ $(0 \leq j \leq d)$ such
that $\alpha_i=0$ for $0 \leq i \leq j-1$ and 
$\alpha_i\not=0$ for $j\leq i \leq d$.
Define
\[
   W = \text{Span}\{\tau_j(A)u, \tau_{j+1}(A)u,\ldots ,\tau_d(A)u\}.
\]
We show that $W$ is nonzero and invariant under each of $A,A^*$.
Observe that $W\not=0$ since $W$ contains $\tau_d(A)u$,
and this vector is nonzero by Lemma \ref{lem:calDu}.
Observe that $AW\subseteq W$, since using (\ref{eq:deftau}) we find
$(A-\theta_iI)\tau_i(A)u = \tau_{i+1}(A)u$ for $j \leq i \leq d-1$ and
$(A-\theta_dI)\tau_d(A)u = 0$.
Observe that $A^*W\subseteq W$, since by (\ref{eq:aux2}) the product
$(A^*-\theta^*_iI)\tau_i(A)u$ is $0$ for $i=j$ and a scalar multiple
of $\tau_{i-1}(A)u$ for $j+1\leq i \leq d$.
We have now shown that $W$ is nonzero and invariant under each of $A,A^*$.
Therefore $W=V$ in view of Definition \ref{def:TDpair}(iv).
We can now easily show that $A,A^*$ is a Leonard pair.
By construction the dimension of $W$ is at most $d+1$.
Also, using (\ref{eq:sumVi}) the dimension of $V$ is $\sum_{i=0}^d \rho_i$,
so this dimension is at least $d+1$ with equality
if and only if $\rho_i=1$ for $0 \leq i \leq d$.
By these comments and since $W=V$ we find
$V$ has dimension $d+1$ and $\rho_i=1$ for $0 \leq i \leq d$. 
Therefore the pair $A,A^*$ is a Leonard pair.

(ii)$\Rightarrow$(i):
Apply Theorem \ref{thm:S} with $X=S$, where $S$ is the
switching element of $A,A^*$.
\end{proof}

\bigskip

\bibliographystyle{plain}

\bigskip\bigskip\noindent
Kazumasa Nomura\\
College of Liberal Arts and Sciences\\
Tokyo Medical and Dental University\\
Kohnodai, Ichikawa, 272-0827 Japan\\
email: knomura@pop11.odn.ne.jp

\bigskip\noindent
Paul Terwilliger\\
Department of Mathematics\\
University of Wisconsin\\
480 Lincoln Drive\\ 
Madison, Wisconsin, 53706 USA\\
email: terwilli@math.wisc.edu

\bigskip\noindent
{\bf Keywords.}
Leonard pair, tridiagonal pair, $q$-Racah polynomial, orthogonal polynomial.

\noindent
{\bf 2000 Mathematics Subject Classification}.
05E35, 05E30, 33C45, 33D45.

\end{document}